\documentclass[a4paper,14pt]{article}

\usepackage{latexsym}
\usepackage{amsmath}
\usepackage{amsfonts}
\usepackage{amssymb}
\usepackage{cite}
\usepackage{mathrsfs}
\newtheorem{lem}{Lemma}[section]
\newtheorem{pro}[lem]{Proposition}
\newtheorem{thm}[lem]{Theorem}

\newtheorem{cor}[lem]{Corollary}

\newcommand{\eps}{\varepsilon}

\renewcommand{\d}{\delta }
\newcommand{\D }{\Delta }

\newcommand{\e }{\varepsilon }

\newcommand{\n }{\nabla }
\newcommand{\vp }{\varphi }

\newcommand{\s }{\sigma }

\renewcommand{\phi}{\varphi}

\renewcommand{\O }{\Omega }
\newcommand{\ov}{\overline}

\newcommand{\be}{\begin{equation}}
\newcommand{\ee}{\end{equation}}
\newenvironment{pf}{\noindent{\sc Proof}.\enspace}{\rule{2mm}{2mm}\medskip}

\newcommand{\R}{\mathbb{R}}
\newcommand{\N}{\mathbb{N}}

\newcommand{\de}{\partial}

\newcommand{\ti}{\tilde}

\newcommand{\M}{\mathcal{M}}
\newcommand{\cM}{\mathcal{M}}

\newcommand{\ra}{{\rangle}}
\newcommand{\la}{{\langle}}

\newcommand{\tn}{\tilde \n}
\DeclareMathOperator{\dist}{dist}
\DeclareMathOperator{\cA}{\mathcal{A}}
\DeclareMathOperator{\diam}{diam}

\title{Profile expansion for the first nontrivial Steklov eigenvalue in Riemannian manifolds}
\author{  Mouhamed Moustapha Fall  and Tobias Weth}

\begin{document}
\date{}
\maketitle
 \let\thefootnote\relax\footnotetext{mouhamed.m.fall@aims-senegal.org (M. M. Fall), weth@math.uni-frankfurt (T. Weth).}

\let\thefootnote\relax\footnotetext{ African Institute for Mathematical Sciences of Senegal,
KM 2, Route de Joal, B. P. 14 18. Mbour, S\'en\'egal.}
\let\thefootnote\relax\footnotetext{ Goethe-Universit\"{a}t Frankfurt, Institut f\"{u}r Mathematik.
Robert-Mayer-Str. 10 D-60054 Frankfurt, Germany.}
\bigskip

\begin{abstract}
We study the problem of maximizing the first nontrivial Steklov eigenvalue of the Laplace-Beltrami Operator among subdomains of fixed volume of a Riemannian manifold. More precisely, we study the expansion of the corresponding profile of this isoperimetric (or isochoric) problem as the volume tends to zero. The main difficulty encountered in our study is the lack of existence results for maximizing domains and the possible degeneracy of the first nontrivial Steklov eigenvalue, which makes it difficult to tackle the problem with domain variation techniques. As a corollary of our results, we deduce local comparison principles for the profile in terms of the scalar curvature on $\cM$. In the case where the underlying manifold is a closed surface, we obtain a global expansion and thus a global comparison principle.\\
\bigskip

\noindent
 \textbf{MSC 2000}: 15A42 - 35B05 - 35R45 - 52A40.
\end{abstract}

\section{Introduction}
Let  $(\M,g)$ be a complete Riemannian manifold of dimension $N\geq2$, and let $\D_g f=div_g(\n f)$ denote the Laplace-Beltrami operator on $\M$. For a bounded regular domain $\O\subset \M$ with outer unit normal $\eta$ on $\de\O$, we consider the Steklov  eigenvalue problem
\begin{equation}
  \label{eq:16}
\D_g f=0\quad \text{in $\O$},\qquad  \la \n f,\eta\ra_g=\nu f\quad
\text{on $\de\O$}.
\end{equation}
The corresponding set of eigenvalues, counted with multiplicities, is given as an increasing sequence
$$
0=\nu_1(\O,g)<\nu_2(\O,g)\leq\dots+\infty.
$$
In the case where $\M= \R^N$, endowed with the euclidean metric $g_{ \text{\tiny $eucl$}}$, it has been proved by Brock \cite{Brock} that, among domains $\O$ of fixed volume $v>0$, balls with volume $v$ are the unique maximizers of $\nu_2(\O)= \nu_2(\O,g_{\text{\rm \tiny $eucl$}})$. In the planar case within the class of simply connected subdomains of $\R^2$, this result had been derived earlier by Weinstock\cite{Weinstock}. The result also extends to the class of simply connected subdomains of a complete Riemannian surface with constant scalar curvature, see \cite[Theorem 7]{Escobar:99}. We point out that, in the euclidean case, Brock \cite{Brock} actually proved the stronger inequality
\begin{equation}
  \label{eq:brock-inequality}
\sum_{i=2}^{N+1} \frac{1}{\nu_i(\O)} \ge N \qquad \text{for every domain $\O$ having the same volume as the unit ball $B \subset \R^N$,}
\end{equation}
with equality if and only if $\O = B$. Note that $\nu_2(B)=\mu_3(B)= \dots= \mu_{N+1}(B)=1$, and the corresponding eigenfunctions on the unit ball are simply the coordinate functions $x \mapsto x^i$, $i=1,\dots,N$. Xia and Wang (see \cite[Theorem 2.1]{Xia-Wang}) also proved a related lower bound for $\sum \limits_{i=2}^{N+1} \frac{1}{\nu_i(\O,g)}$ in the case where $(\cM,g)$ is a Hadamard manifold. \\
In the present paper we  study the geometric variational problem of maximizing
$\nu_2(\O,g)$ among domains with fixed small volume in a
general Riemannian manifold $(\M,g)$.  For $0<v< |\M|_g$, we define the
\textit{Weinstock-Brock profile} of $\M$ as
$$
{W\!B}_{\M}(v,g):=\sup_{\O\subset\M,\,|\O|_g=v}\,\nu_2(\O,g).
$$
Here and in the following, we assume without further mention that only
regular bounded domains $\O \subset \M$ are considered, and we let $|\Omega|_g$ denote the $N$-dimensional volume with respect to the metric $g$. For open subsets $\cA
\subset \M$ and $0<v< |\cA|_g$, we also define
$$
{W\!B}_{\cA}(v,g):=\sup_{\O\subset\cA,\,|\O|_g=v}\,\nu_2(\O,g),
$$
assuming again without further mention that only regular bounded
domains $\O \subset \cA$ are considered. By Brock's result \cite{Brock} mentioned above and the scaling properties of $\nu_2$, we then have
$$
{W\!B}_{\R^{N}}(v)=\left(\frac{v}{|B|}\right)^{-\frac{1}{N}}.
$$
In our first result we analyze the local
effect of the scalar curvature of $\M$ on the $\nu_2$-profile. For this
we let $B_g(y_0,r)$ denote the geodesic ball in $\M$
centered at a point $y_0 \in \M$ with radius $r$. The following result contains a global asymptotic lower bound for ${W\!B}_{\M}(v)$ and a sharp two-sided bound for
${W\!B}_{B_g(y_0,r)}(v)$ if $r>0$ is small.

\begin{thm}\label{pro:Main-rslt}
Let $\M$ be a complete $N$-dimensional Riemannian manifold with
$N\geq 2$, and let $S$ denote the scalar curvature function on $\M$.
Moreover, let $y_0 \in \M$. Then we have:
\begin{itemize}
\item[(i)] As $v \to 0$,
\begin{equation}
  \label{eq:24}
{W\!B}_{\M}(v) \;\ge \;\Bigl(\frac{v}{|B|}\Bigr)^{-\frac{1}{N}} +\frac{S(y_0) }{2N(N+2) }
\Bigl( \frac{v}{|B|} \Bigr)^{\frac{1}{N}}
+o(v^{\frac{1}{N}}).
\end{equation}
\item[(ii)] For every $y_0 \in \M$ and every $\e>0$ , there exists $r_\e>0$ such that
\begin{equation}
\!\!\!\!\!\!\!\!{W\!B}_{\!B_g(y_0,r_\e)}(v)  \: \left\{
  \begin{aligned}
 &\!\!\ge \Bigl(\frac{v}{|B|}\Bigr)^{-\frac{1}{N}}+\left(\frac{S(y_0) }{2N(N+2) } -\e
\right) \Bigl( \frac{v}{|B|} \Bigr)^{\frac{1}{N}}\\
 &\!\!\leq\Bigl(\frac{v}{|B|}\Bigr)^{-\frac{1}{N}}+ \left(\frac{S(y_0) }{2N(N+2) } +\e \right)
 \Bigl(\frac{v}{|B|} \Bigr)^{\frac{1}{N}}
  \end{aligned}
\right. \quad \text{for $v\in \bigl(0\,,\,| B_g(y_0,r_\e)|_g\bigr)$.}
\end{equation}
\end{itemize}
\end{thm}
We note that $S(y_0)$ in (\ref{eq:24}) can be replaced by $\sup_{\M}S$ if the supremum is attained on $\M$ (e.g. if $\M$ is compact). The result naturally leads to the question whether a sharp upper bound can also be obtained for ${W\!B}_{\M}(v)$. The main problem which arises here is the fact that almost maximizing domains of small volume $v$ do not necessarily have small diameter if $N \ge 3$.  However, we are able to control the diameter in the two-dimensional case, and thus we have the following result.

\begin{thm}
\label{sec:global-upper-bound-5-0}
Let $(\M,g)$ be a closed Riemannian surface.  Then we have
$$
{W\!B}_{\M}(v) \;= \;\Bigl(\frac{v}{\pi}\Bigr)^{-\frac{1}{2}} +\frac{S_{\M}}{16} \Bigl(\frac{v}{\pi}\Bigr)^{\frac{1}{2}}
+o(v^{\frac{1}{2}})\qquad \text{as $v \to 0$, }
$$
where $S_{\M}$ denotes the maximum of the scalar curvature function $S$ on $\M$.
\end{thm}

We conjecture that a similar global expansion holds in closed Riemannian manifolds of higher dimension, but for now this remains open. As we will explain below in more detail, our proof of Theorem~\ref{sec:global-upper-bound-5-0} does not extend to higher dimensions. 

 An immediate consequence of the asymptotic estimates given in Theorem~\ref{pro:Main-rslt} and \ref{sec:global-upper-bound-5-0} are the following comparison principles.

\begin{cor}\label{th:main-th-2d}
Let $(\M_1,g_1)$, $(\M_2,g_2)$ be two $N$-dimensional complete
Riemannian manifolds, $N \ge 2$ with scalar curvature functions
$S_1$, $S_2$ respectively.
\begin{itemize}
\item[(i)] Let $y_1\in\M_1$ and $y_2 \in \M_2$ such
that $S_1(y_1)<S_2(y_2)$. Then there exists $r>0$ such that
  \begin{equation*}
{W\!B}_{B_{g_1}(y_1,r)}(v) < {W\!B}_{B_{g_2}(y_2,r)}(v)\qquad \text{for any $v \in
(0,\min\{|B_{g_1}(y_1,r)|_{g_1},|B_{g_2}(y_2,r)|_{g_2}\}).$}
  \end{equation*}
\item[(ii)] If $N=2$ and $(\M_1,g_1)$, $(\M_2,g_2)$ are closed Riemannian surfaces with $\max \limits_{M_1}S_1 < \max \limits_{M_2} S_2$, then there exists $r>0$ such that
\begin{equation*}
{W\!B}_{\M_1}(v) < {W\!B}_{\M_2}(v)\qquad \text{for any $0 < v < r$.}
  \end{equation*}
\end{itemize}
\end{cor}

Our results should be seen in comparison with our recent work \cite{FW-Neum} on the Szeg\"o-Weinberger profile in Riemannian manifolds, which arises from the corresponding maximization problem for the first nontrivial Neumann eigenvalue of $-\Delta_g$ on $\cM$. In this work we established an analogue of Theorem~\ref{pro:Main-rslt} for the the Szeg\"o-Weinberger profile. Similarly as in \cite{FW-Neum}, the first step in the proof of Theorem~\ref{pro:Main-rslt} is the derivation of expansions for $\nu_2$ for small geodesic balls and small ellipsoids with small eccentricity centered at a point $y_0 \in \cM$. More precisely, in Proposition~\ref{lem:expa_un2B} below we show that
\begin{equation}
  \label{eq:exp-ball-0}
\nu_2(B_g(y_0,r),g) =\frac{1}{r}+
\frac{2r}{3(N+2)}R_{min}(y_0)+o(r) \qquad \text{as $r \to 0$}
\end{equation}
with $R_{min}(y_0):=\min \limits_{A\in T_{y_0}\M, |A|_g=1}Ric_{y_0}(A,A)$. Hence there is an anisotropic curvature effect on the expansion which suggests that small geodesic balls are not optimal up to linear order in $r$ for the maximization problem. We therefore construct a family (depending on $r$) of small ellipsoids $E(y_0,r)$ which are choosen such that the eccentricity balances the anisotropic curvature effects, so that the resulting expansion
\begin{equation}
  \label{eq:29-0}
\nu_2(E(y_0,r),g) =\frac{1}{r}+\frac{2r}{3N(N+2) }S(y_0)+o(r)
\end{equation}
depends only on the scalar curvature $S(y_0)$, see Proposition~\ref{sec:expansion-mu_2-small-2} below.
The computations of these expansions bear some similarities with the corresponding ones in \cite{FW-Neum}, although some differences arise due to the fact that boundary integrals have to be expanded in the present case. On the other hand, we note that the simple form of the eigenfunctions corresponding to $\nu_2(B)$ leads to a nicer expansion than in the Neumann eigenvalue case. We shall see that by combining (\ref{eq:29-0}) with the volume expansion for $E(y_0,r)$, we already obtain the lower bound for the profile given in Theorem~\ref{pro:Main-rslt}(i).  The proof of the local upper bound in Theorem~\ref{pro:Main-rslt}(ii) is more involved and proceeds eventually by a contradiction argument. For this, some care is needed to construct, for given subdomains of $B_g(y_0,r)$ with $r>0$ small, suitable vector fields which can be used in combination with the variational principle for $\nu_2$ in order to control the symmetric distance of these domains to a suitably chosen geodesic ball with the same volume. Within this step, the key tool is a quantitative weighted isoperimetric inequality proved recently by Brasco, de Philipps and Ruffini   see \cite[Theorem B]{BDR}.\\
 We point out that, in the proof of the local upper bound for the profile given in Theorem~\ref{pro:Main-rslt}(ii),  the arguments differ significantly from the ones in \cite{FW-Neum} for the Neumann eigenvalue case. We also remark that, at least up to now, Theorem~\ref{sec:global-upper-bound-5-0} has no analogue for the corresponding Neumann eigenvalue profile. The proof of this global expansion is technically involved, but the strategy is easy to explain. We will show that almost maximizing domains for $\nu_2$ of small (fixed) volume must also have small diameter. There is no hope to prove this in dimension $N \ge 3$, since in this case one may increase the diameter of the domain by adding a long cusp of small volume and perimeter. By the variational characterization, this will only result in a small change of $\nu_2$. In contrast, as remarked before, in the two-dimensional case we will be able to deduce bounds on the diameter with the help of the variational characterization of $\nu_2$ and suitably constructed test functions.\\
To close the introduction, we mention the earlier work in \cite{Druet-FK,Fall-eigen} on the small volume expansion for the Faber-Krahn profile, which is related to the minimization of the first Dirichlet eigenvalue $\lambda_1(\Omega,g)$ of $-\Delta_g$ among subdomains $\Omega$ of fixed volume. One important difference between $\lambda_1(\Omega,g)$ and $\nu_2(\Omega,g)$ is the degeneracy of $\nu_2$ in the case of the unit ball and possibly also in the case of maximizing domains on Riemannian manifolds. This degeneracy makes it difficult to apply domain variation arguments to the maximization problem.

The paper is organized as follows. Section \ref{s:pn} contains some preliminaries and the proof of local expansions of $\nu_2$ for small geodesic balls and small ellipsoids with small eccentricity. In particular, as already remarked above, we shall see that suitably chosen ellipsoids provide the optimal lower bound in Theorem~\ref{pro:Main-rslt}(i). In Section~\ref{s:Ubm2} we then complete the proof of Theorem~\ref{pro:Main-rslt} by providing the upper bound in (ii). Finally, in Section~\ref{sec:global-upper-bound}, we focus on the two-dimensional case $N=2$ and give the proof Theorem~\ref{sec:global-upper-bound-5-0}.\\

\noindent

\textbf{Acknowledgments:}
The authors would like to thank Lorenzo Brasco for  many helpful comments on previous 
versions of the paper, stimulating discussions, and for drawing our attention to the paper \cite{BDR}. The first author  is supported by the Alexander von Humboldt foundation. 
Part of the paper was written while the second author was visiting AIMS Senegal in November 2014 
within the project {\em ``Joint steps in geometric variational problems''} funded by DAAD, Germany. He wishes to thank AIMS Senegal for the kind hospitality and DAAD for the funding of the visit.
\\

\noindent {\em General Notation:}
Throughout the paper, we let  $B$ denote the open unit ball in $\R^{N}$ and put
$rB:= \{x \in \R^N\::\: |x|<r\}$ for $r>0$. Moreover, we write $x \cdot
y$ for the euclidean scalar product of $x,y \in \R^N$.

\section{Local expansions of $\nu_2$ for small geodesic balls and ellipsoids.}\label{s:pn}
Let  $(\M,g)$ be a complete Riemannian manifold of dimension $N \ge 2$. For a smooth
bounded subdomain domain $\Omega$ of $(\M,g)$,
we write $\nu_2= \nu_2(\Omega,g)$ for the first nontrivial
eigenvalue of (\ref{eq:16}). The variational characterization of
$\nu_2(\Omega,g)$ is given by
\begin{equation}
  \label{eq:46}
\nu_2(\Omega,g) = \inf \Bigl\{\int_{\Omega}|\nabla u|_g^2d
v_g\::\: u \in H^1(\Omega),\:\int_{\de \Omega}
  u^2 d\sigma_g=1,\:
\int_{\de \Omega}u\,d \sigma_g=0 \Bigr\}.
\end{equation}
Here $v_g$ denotes the volume element of the metric $g$, and $\sigma_g$ denotes the volume element of the restriction of $g$ to an $N-1$-dimensional submanifold of $\M$. For a Borel subset $A
\subset \M$, we let $|A|_g$ denote the $N$-dimensional volume of $\Omega$ and $\sigma_g(A)$ denote the $N-1$-dimensional Haussdorff-measure, both with respect to the metric $g$.
 If $\M= \R^{N}$
and $g$ is the Euclidean metric, we simply write $dx$ in place of $dv_g$, $|\cdot|$ in place of $|\cdot|_g$, $d \sigma$ in place of $d \sigma_g$ and $\nu_2(\Omega)$ in
place of $\nu_2(\Omega,g)$. We recall
that the minimizers of the minimization problem (\ref{eq:46}) are precisely the
eigenfunctions corresponding to $\nu_2(\Omega,g)$.  As noted already, in the case of the unit ball $B \subset \R^N$ we have that $\nu_2(B)=1$ is of
multiplicity $N$ with corresponding eigenfunctions given by $x \mapsto x^i$, $i=1,\dots,N$.\\

In the following, we assume that $(\M,g)$ is complete, and we
fix $y_0\in\M$ and an orthonormal basis $E_1,\dots,E_{N}$
of $T_{y_0}\M$.  We will use the
(somewhat sloppy) notation
$$
X:=x^i E_i \in T_{y_0}\M \qquad \text{for $x \in \R^{N}$.}
$$
Here and in the following, we sum over repeated upper and lower indices as usual. We
consider the map
\be\label{eq:sys-coor}
\Psi: \R^{N} \to \M, \qquad  \Psi(x):=\textrm{Exp}_{y_0}(X),
\ee
which gives rise to a local geodesic coordinate system of a neighborhood of $y_0$. A geodesic ball in $\M$ centered at $y_0$ with radius $r>0$ is given
as $B_g(y_0,r)=\Psi(rB)$. The map $\Psi$ induces coordinate
vector fields $Y_i:=\Psi_*  \frac{\partial}{\partial x^i}$ on $\R^N$ given by
$$
Y_i(x)= d\, \textrm{Exp}_{y_0}(X) E_i \in T_{\Psi(x)} \M, \qquad \text{for $x \in \R^N$, $i=1,\dots,N.$}
$$
We need local expansions for the associated  metric coefficients
$$
g_{ij}(x)=\la Y_i(x),Y_j(x)\ra_g \qquad \text{for $x \in \R^{N}$, $i.j =1,\dots,N$.}
$$
For this we let $R_{y_0}: T_{y_0}\M \times T_{y_0}\M \times T_{y_0}\M \to T_{y_0}\M$
denote the Riemannian curvature tensor at $y_0$ and
$$
Ric_{y_0} : T_{y_0}\M \times T_{y_0}\M \to \R,\qquad Ric_{y_0}(X,Y)=-\sum_{i=1}^{N}\la R_{y_0}(X,E_i)Y,E_i\ra_{g}
$$
the Ricci tensor at $y_0$. Moreover, we let $S:\M \to \R$ denote the scalar curvature function
on $\M$, so that $S(y_0)= \sum \limits_{k=1}^N Ric_{y_0}(E_k,E_k)$. It will be useful to put
\begin{equation}
  \label{eq:17}
R_{ijkl}:= \la R_{y_0}(E_i,E_j)E_k,E_l\ra_g  \quad \text{and}\quad R_{ij}:= Ric_{y_0}(E_i,E_j) \qquad \text{for $i,j=1,\dots,N$.}
\end{equation}
Without changing the value of these constants,  we sometimes raise lower to upper indices in the following.  We then have the  following well known local expansions as $|x| \to 0$ (see e.g.
in \cite[\S II.8]{chavel-1}):
\begin{align}
 g_{ij}(x)&=\delta_{ij}+\frac{1}{3}\,\la R_{y_0}(X,E_i)X,E_j\ra_g
+ {O}(|x|^3) = \delta_{ij}+\frac{1}{3}R_{kilj}x^kx^l
+ {O}(|x|^3); \label{exp-fund-1}\\
dv_{g}(x)&=\Bigl( 1-\frac{1}{6}\, Ric_{y_0} (X,X) + O(|x|^3) \Bigr)dx=\Bigl( 1-\frac{1}{6}R_{lk}x^l x^k + O(|x|^3) \Bigr)dx.\label{exp-fund-2}
\end{align}
As a consequence of (\ref{exp-fund-2}), the volume
expansion of metric balls is given by
\be\label{eq:expvolBgr}
\left|B_g(y_0,r) \right|_g= r^{N}\,{\left|B \right|}\,\left(1-\frac{1}{6(N+2)}\,r^2{S}(y_0)+O(r^4)  \right).
\ee
The first goal of this section is to derive an expansion for $\nu_2$ on small geodesic balls centered at
$y_0$.  It will be useful to pull back the problem to the unit ball $B \subset \R^{N}$.
For this we let $r>0$ be smaller than half of the injectivity radius of $\M$ at $y_0$, so
that $B_g(y_0,s)$ is a regular domain for $s \le 2r$.  Moreover, we consider the pull back metric
of $g$ under the map $2B \to \M,\;x \mapsto \Psi(rx)$, rescaled with the
factor $\frac{1}{r^2}.$ Denoting this metric on $2B$ by $g_r$, we then
have, in euclidean coordinates,
$$
[g_r]_{ij}(x)= \langle \frac{\de}{\de x^i},\frac{\de}{\de x^j}
\rangle_{g_r}\Big|_x  = \langle Y_i(\Psi(rx)),
Y_j(\Psi(rx))\rangle_{g}= g_{ij}(rx),
$$
so that, as a consequence of (\ref{exp-fund-1}),
\begin{equation}
  \label{eq:1}
[g_r]_{ij}(x)=\d_{ij}+\frac{r^2}{3}R_{kilj}x^k x^l+O(r^3) \qquad \text{as $r \to 0$}
\end{equation}
and
\begin{equation}
\label{eq:2}
g_r^{ij}(x)=\d^{ij}-\frac{r^2}{3}R_k\!{}^i\!{}_l{}^jx^k x^l+O(r^3) \qquad \text{as $r \to 0$}
\end{equation}
uniformly for $x \in \overline B$. Here, as usual,
$(g_r^{ij})_{ij}$ denotes the inverse of the matrix
$([g_r]_{ij})_{ij}$. Setting $|g_r|= \det
([g_r]_{ij})_{ij}$, we also have
\begin{equation}
  \label{exp|g_r|}
\sqrt{|g_r|}(x) = 1-\frac{r^2}{6}
R_{kl} x^kx^l +O(r^3) \qquad \text{as $r \to 0$}
\end{equation}
uniformly for $x \in \overline B$ by (\ref{exp-fund-2}). Since this expansion is valid
in the sense of $C^1$-functions on $\overline B$, it follows that
\begin{equation}
\frac{\de}{\de x^i}\sqrt{|g_r|}= -\frac{r^2}{3}R_{ki}x^k+O(r^3) \qquad \text{as $r \to 0$ for $i=1,\dots,N$.}  \label{eq:3}
\end{equation}
The expansion~(\ref{exp|g_r|}) obviously yields
\begin{equation}
  \label{exp|g_r|-volume}
dv_{g_r}(x)=\Bigl(1-\frac{r^2}{6}R_{lk}x^l x^k + O(r^3) \Bigr)dx\qquad \text{as $r \to 0$} 
\end{equation}
uniformly for $x \in \overline B$. We will also need the following expansion for boundary integrals with respect to subdomains of $B$.

\begin{lem}
\label{sec:local-expans-nu_2-surface}
For every smooth domain $U \subset B$ and every $f\in C^1(\de U)$ we have
\begin{equation}
  \label{eq:boundary-integral}
\int_{\de U}f(x)  \,d\sigma_{g_{r}}= (1+O(r^2)) \int_{\de U}f(x) d\s,
\end{equation}
where $\frac{O(r^2)}{r^2}$ remains bounded uniformly in $U$ and $f$ as $r \to 0$. Moreover, for every $f\in
C^1(\de B)$ we have
\begin{equation}
  \label{exp|g_r|-surface}
\int_{\de B}f(x)  \, d\sigma_{g_r}(x)= \int_{\de B} \Bigl(1-\frac{r^2}{6}R_{lk}x^l x^k\Bigr)f(x)d\s  + O(r^3)\int_{\de B}f(x) d\s,
\end{equation}
where $\frac{O(r^3)}{r^3}$ remains bounded uniformly in $f$ as $r \to 0$.
\end{lem}

We note that (\ref{exp|g_r|-surface}) follows from the computations in \cite[Appendix 4.1]{pacard-xu}. Here we provide a different short proof, based on integration by parts.\\

\begin{pf}
Let $\eta_r$ denote the unit outer normal vector field on $\partial U$ with respect to $g_r$ and $\eta$ the unit outer normal vector field on $\partial U$ with respect to the euclidean metric. We first claim that, for fixed $r>0$,
\be \label{eq:intU_rdsrf}
\int_{\de U }f  \,d\sigma_{g_{r}}=  \int_{\de U_{r}} f \sqrt{|g_{r}|} \, \eta_r \!\cdot \!\eta  \, d\s \quad \textrm{ for every } f\in C^1(\partial U),
\ee
where, as before, $\cdot$ denotes the euclidean scalar product. To show this, we may first extend $f \eta_r: \partial U \to \R^N$ to a $C^1$-vector field $\xi$ on $\R^N$.
Applying the divergence theorem with respect to the metric $g_r$, we then have
\begin{equation}
\int_{U}\textrm{div}_{g_r}\xi \:\sqrt{|g_{r}|}dx=\int_{U}\textrm{div}_{g_r}\xi \:dv_{g_r}= \int_{\de U} \la \xi, \eta_r \ra_{g_r}  \,d\sigma_{g_{r}}=
 \int_{\de U} f  \,d\sigma_{g_{r}}.
\end{equation}
On the other hand, applying the divergence theorem with respect to the euclidean metric, we find that
$$
\int_{U}\textrm{div}_{g_r}\xi \:\sqrt{|g_{r}|}dx=
\int_{U} \frac{\de}{\de x^i}\Bigl[\xi^i \sqrt{|g_{r_k}|}\Bigr] dx= \int_{\de U} \sqrt{|g_{r_k}|} \:\xi \!\cdot\! \eta\: d\s =
\int_{\de U} f \sqrt{|g_{r_k}|} \: \eta_r \!\cdot\! \eta\: d\s.
$$
Hence \eqref{eq:intU_rdsrf} follows. In order to expand the term $\eta_r \!\cdot\! \eta$ in $r$, we consider a point $q\in\de U$ and let $e_i$, $i=1,\dots, N-1$ be an orthonormal basis of $T_q\de U$ with respect to the Euclidean metric. For simplicity, we will write $\eta_r$ and $\eta$ instead of $\eta_r(q)$ and $\eta(q)$ in the following. We then have that
\begin{equation}
  \label{eq:additional-boundary-orthonormal}
\eta_r=[\eta_r \!\cdot\! \eta] \eta  + [\eta_r  \!\cdot\! e_i] e_i
\end{equation}
Since $e_i\in T_q \de U$, we have, by (\ref{eq:1}),
\begin{equation}
  \label{eq:additional-boundary-orthonormal-(1)}
0 = \la  \eta_r,e_i \ra_{g_{r}}= \eta_r  \!\cdot\! e_i +O(r^2)|\eta_r| |e_i| = \eta_r  \!\cdot\! e_i + O(r^2)|\eta_r|.
\end{equation}
Moreover, (\ref{eq:1}) also implies that
\begin{equation}
  \label{eq:additional-boundary-orthonormal-0}
1= |\eta_r|_{g_r}^2 = |\eta_r|^2 + \frac{r^2}{3} R_{kilj}\, \eta_r^i \eta_r^j q^k q^l  +O(r^3)|\eta_r|^2  = (1+O(r^2))|\eta_r|^2
\end{equation}
and hence
\begin{equation}
  \label{eq:additional-boundary-orthonormal-1}
|\eta_r| = 1+O(r^2) 
\end{equation}
Consequently, (\ref{eq:additional-boundary-orthonormal-(1)}) implies that $\eta_r \cdot e_i= O(r^2)$ independently of $U \subset B$, $q$ and the choice of the orthonormal basis $e_i$. Taking the euclidean scalar product of (\ref{eq:additional-boundary-orthonormal}) with $\eta_r$, we now find that
\begin{equation}
  \label{eq:additional-boundary-orthonormal-3}
|\eta_r|^2 = [\eta_r \!\cdot\!  \eta ]^2 +  O(r^4).
\end{equation}
Together with (\ref{eq:additional-boundary-orthonormal-1}) this implies that $[\eta_r \!\cdot\!  \eta ]^2 = 1 + O(r^2)$.
Since both $\eta_r$ and $\eta$ are defined as outer normal vector fields, we conclude that
\begin{equation}
  \label{eq:additional-boundary-orthonormal-2}
\eta_r \!\cdot\!  \eta = 1 + O(r^2) \qquad \text{uniformly on $\partial U$ and independently of $U$}.
\end{equation}
Combining this with (\ref{exp|g_r|}) and \eqref{eq:intU_rdsrf}, we obtain (\ref{eq:boundary-integral}).\\
To see (\ref{exp|g_r|-surface}), we consider the special case $U=B$, and we note that, as a consequence of Gauss' Lemma (see e.g. \cite[Corollary 5.2.3]{jost:RiemannianGeometryandGeometricAnalysis}), we have that $\eta_r(q)= q= \eta(q)$ for every $q \in \partial B$. Together with (\ref{exp|g_r|}) and \eqref{eq:intU_rdsrf} this implies~(\ref{exp|g_r|-surface}).
 \end{pf}

 We are now ready to establish the desired expansion for $\nu_2$ on small geodesic balls.

\begin{pro}\label{lem:expa_un2B}
We have
\begin{equation}
  \label{eq:exp-ball}
\nu_2(B_g(y_0,r),g) =\frac{1}{r}+
\frac{2r}{3(N+2)}R_{min}(y_0)+o(r) \qquad \text{as $r \to 0$}
\end{equation}
with $R_{min}(y_0):=\min \limits_{A\in T_{y_0}\M, |A|_g=1}Ric_{y_0}(A,A).$
\end{pro}
\begin{pf}
Let $u_r \in C^3(\overline{B_g(y_0,r)})$ be an eigenfunction associated to $\nu_2(B_g(y_0,r),g)$, i.e., we have
$$
\D_g u_r =0\quad \text{in $B_g(y_0,r)$},\qquad  \la \n
u_r,\eta_r\ra_g=\nu_2(B_g(y_0,r),g) \,u_r\quad \text{on $\de
B_g(y_0,r)$,}
$$
where $\eta_r$ denotes the outer unit normal on $\de B_g(y_0,r)$.  Then the function
$\Phi_r: \overline B \to \R$, $\Phi_r(x)=r^{\frac{N}{2}}u_r(\Psi(rx))$ satisfies
 \be\label{eq:DPhir}
  \D_{g_r} \Phi_r=0\quad \text{in $B$},\qquad \la
\n \Phi_r,\eta \ra_{g_r}=\nu_2(B,g_r ) \,\Phi_r\quad \text{on $\de
B$,}
 \ee
  with
$$
\D_{g_r}\Phi_r=  \frac{1}{{\sqrt{|g_r|}}}\frac{\de }{\de
x^i}\left(\sqrt{|g_r|}g_r^{ij} \frac{\de \Phi_r}{\de x^j}\right)
\qquad \text{and}\qquad \nu_2(B,g_r)=r \nu_2(B_g(y_0,r),g).
$$
Hence the asserted expansion (\ref{eq:exp-ball}) is equivalent to
\begin{equation}
  \label{exp-equiv}
 \nu_2(B,g_r)= 1+
\frac{2r^2}{3(N+2)}R_{min}(y_0)  +o(r^2).
\end{equation}
Moreover, by normalization we may assume that  $\int_{\de B} \Phi_r^2\,d\sigma_{g_r}=1$. To prove (\ref{exp-equiv}), we first note that, since
$g_r$ converges to the Euclidean metric in $\ov{B}$, it follows
from the variational characterization of $\nu_2$ that
$\nu_2(B,g_r )\to \nu_2(B)=1$. Moreover, by using standard elliptic
regularity theory, one may show that,
along a sequence $r_k \to 0$, we have
$\Phi_{r_k}\to \Phi$ in $H^1(B)$ for some function $\Phi \in
C^2_{loc}(B)\cap C^1(\ov{B})$ satisfying
\begin{equation}
  \label{limiteq-1}
\D \Phi=0\quad \text{in $B$} ,\qquad  \la \n \Phi,\eta
  \ra= \Phi\quad \text{on $\de B$},  \quad
  \text{and}\quad \int_{\partial B} \Phi^2\,d\sigma=1.
\end{equation}
Hence there exists $a=(a_1,\dots,a_N)=(a^1,\dots,a^N) \in \R^{N}$
with $|a|=1$ and such that
\begin{equation}
  \label{eq:def-phi}
\Phi(x)= \frac{a  \cdot x }{\sqrt{|B|}} \qquad \text{for $x  \in \overline B$.}
\end{equation}
For matters of convenience, we will continue to write $r$ instead of
$r_k$ in the following. By integration by parts in \eqref{eq:DPhir},
using $\la \n\Phi_r,\eta\ra_{g_r}=0$ and $dv_{g_r}= \sqrt{|g_r|}dx$,
we have
$$
\nu_2(B, g_r)\, \int_{\de B} \Phi \Phi_r d\sigma_{g_r}=
 \int_{B}\sqrt{|g_r|}g_r^{ij} \frac{\de \Phi}{\de x^i}
                                          \frac{\de \Phi_r}{\de x^j} \,dx.
$$
We thus find, using (\ref{eq:2}), (\ref{eq:3}), (\ref{limiteq-1}) and integrating by parts again,
\begin{align}
&\int_{B}\sqrt{|g_r|}g_r^{ij} \frac{\de \Phi}{\de x^i}
                                           \frac{\de \Phi_r}{\de x^j}
                                           \,dx
= \int_{B} \sqrt{|g_r|}\Bigl( \nabla \Phi_r \cdot \nabla \Phi -\frac{r^2}{3}   R_k\!{}^i\!{}_l{}^jx^kx^l  \frac{\de \Phi}{\de x^i}\frac{\de \Phi_r}{\de x^j}
\Bigr)dx+O(r^3)
\label{eq:4}\\
&=\int_{\de B}\Phi\Phi_r\,d\sigma_{g_r}
                                  - \int_{B}\Phi_r  \nabla
                                  \sqrt{|g_r|} \cdot \nabla \Phi  \,dx
                                          -\,\frac{r^2}{3} \int_{B}\la
                                          R_k\!{}^i\!{}_l{}^jx^kx^l  \frac{\de \Phi}{\de x^i}\frac{\de \Phi_r}{\de x^j}\,dx+O(r^3)
                                          \nonumber \\
&=\int_{\de
B}\Phi\Phi_r\,d\sigma_{g_r}+\frac{r^2}{3}\int_{B}\Phi_r\,
R_{i}{}^jx^i \frac{\de \Phi}{\de x^j} \,dx-\,\frac{r^2}{3} \int_{B}R_k\!{}^i\!{}_l{}^jx^kx^l  \frac{\de \Phi}{\de x^i}\frac{\de \Phi_r}{\de x^j}
\,dx+O(r^3).\nonumber
\end{align}
Therefore, since $ \int_{\de B}\Phi\Phi_r\,d\sigma_{g_r}\to  1$ and
$\Phi_r\to \Phi$ in $H^1(B)$ as $r\to0$, we obtain
\begin{align}
 \nonumber
\nu_2(B, g_r)&=1+\frac{r^2}{3}\int_{B}\Phi\,
R_{i}{}^jx^i \frac{\de \Phi}{\de x^j} \,dx-\,\frac{r^2}{3} \int_{B}R_k\!{}^i\!{}_l{}^jx^kx^l  \frac{\de \Phi}{\de x^i}\frac{\de \Phi}{\de x^j}
\,dx+o(r^2)\\
&=
1+\frac{r^2}{3|B|}\int_{B} a_k x^k
 R_{i}{}^jx^i a_j \,dx-\,\frac{r^2}{3|B|} \int_{B} R_k\!{}^i\!{}_l{}^j x^kx^l  a_i a_j\,dx+o(r^2).
 \label{eq:34}
\end{align}
Recalling that
\begin{equation}
  \label{Bxixj-int}
\int_{B}x^i x^j\,dx =\delta^{ij}\frac{|B|}{N+2} \qquad \text{for $i,j= 1,\dots,N$,}
\end{equation}
we calculate
\begin{equation}
  \label{eq:25-new}
\int_{B} a_k x^k
 R_{i}{}^jx^i a_j \,dx   = \frac{|B|}{N+2}  R^{kj} a_k a_j  = \frac{|B|}{N+2}Ric_{y_0}(A,A)
 \end{equation}
with $A:= a_i E^i$ and
\begin{equation}
  \label{eq:27-new}
\int_{B}R_k\!{}^i\!{}_l{}^jx^kx^l  a_i a_j\,dx =  -\frac{|B|}{N+2} R^{ij}  a_i a_j
=  -\frac{|B|}{N+2} Ric_{y_0}(A,A).
\end{equation}
Therefore \begin{equation}
  \label{eq:10}
\nu_2(B,g_r)= 1+
\frac{2r^2}{3(N+2)}Ric_{y_0}(A,A)  +o(r^2).
\end{equation}
We now need to recall that -- more precisely -- here we consider a sequence $r=r_k \to 0$.
Nevertheless, the argument implies that
\begin{equation}
  \label{eq:8}
\nu_2(B,g_r) \ge 1 +
\frac{2r^2}{3(N+2)}R_{min}(y_0)  +o(r^2)  \qquad \text{as $r \to 0$.}
\end{equation}
Indeed, if - arguing by contradiction - there is
a sequence $r_k \to 0$ such that
\begin{equation}
  \label{eq:9}
\limsup_{k \to \infty} \nu_2(B,g_r)  <   1+
\frac{2r_k^2}{3(N+2)}R_{min}(y_0)  +o(r_k^2)  \qquad \text{as $k \to \infty$.}
\end{equation}
then, by the above argument, there exists a subsequence along which the
expansion (\ref{eq:10}) holds with some $A \in T_{y_0}\M$ with $|A|=1$, thus
contradicting (\ref{eq:9}). Hence (\ref{eq:8}) is true, and it thus remains to prove that
\begin{equation}
  \label{eq:44}
\nu_2(B,g_r) \leq 1
+\frac{2r^2}{3(N+2)}Ric_{y_0}(A,A)  +o(r^2
) \quad \text{for all $A \in T_{y_0}\M$ with $|A|=1$.}
\end{equation}
 So consider $a=(a^1,\dots,a^N) \in \R^N$ arbitrary with
$|a|=1$, let $A=a^i E_i \in
T_{y_0}\M$, and  define $\Phi: \overline B \to \R$ by (\ref{eq:def-phi}). Moreover, put
$c_r:= \frac{1}{|\de B|_{g_r}}\int_{\de B} {\Phi}d\sigma_{g_r}$ for $r>0$ small. Then, by (\ref{exp|g_r|-surface}),
$$
c_r = \Bigl(\frac{1}{|\de B|}+O(r^{2})\Bigr)\Bigl(
 \int_{\de B}{\Phi}(x)[1-\frac{1}{6} R_{kl}x^kx^l]d\sigma(x)
 +O(r^3)\Bigr)=\Bigl(\frac{1}{|\de B|}+O(r^{2})\Bigr)O(r^3)=O(r^3),
$$
since the function $x \mapsto {\Phi}(x)[1-\frac{1}{6} R_{kl}x^kx^l]$ is odd with
respect to reflection at the origin. Hence, using the
variational characterization of $\nu_2(B,g_r)$, we find that
$$
\nu_2(B,g_r)\leq\,\frac{  \displaystyle \int_{B } |\n
({\Phi}-c_r)|_{g_r}^2\, d{v_{g_r}}    }{\displaystyle \int_{\de
B}
 ({\Phi}-c_r)^2 \,d{\sigma_{g_r}}}
=\,\frac{  \displaystyle \int_{B } |\n {\Phi}|_{g_r}^2\,
d{v_{g_r}} }{\displaystyle  \int_{\de B} \bigl[{\Phi}^2 +O(r^3)\bigr]
\,d{\sigma_{g_r}}} =\,\frac{  \displaystyle \int_{B } |\n
{\Phi}|_{g_r}^2\,d{v_{g_r}}  +O(r^3)}{\displaystyle  \int_{\de B}
 {\Phi}^2 \,d{\sigma_{g_r}}}
$$
and therefore
$$
\nu_2(B,g_r) \int_{\de B }
  {\Phi}^2\,d{\sigma_{g_r}} \le
 \int_{B_g(y_0,r) } |\n {\Phi}|_{g_r}^2\,d{v_{g_r}} +O(r^3)=
\int_{B}\sqrt{|g_r|}g_r^{ij} \frac{\de \ti {\Phi}}{\de x^i}
                                           \frac{\de {\Phi}}{\de x^j}
                                           \,dx  +O(r^3).
$$
It is by now straightforward that the same estimates as above -- starting from
(\ref{eq:4}) -- hold with both $\Phi_r$ and $\Phi$ replaced by
${\Phi}$. We thus obtain (\ref{eq:44}), as required.  
\end{pf}

\begin{cor}
\label{sec:expansion-mu_2-small-1} We have \be
\label{eq:expmu_2}
\nu_2(B_g(y_0,r),g)=\Bigl(\frac{v}{|B|}\Bigr)^{-\frac{1}{N}} + \frac{4NR_{min}(y_0)-S(y_0)}{ 6N(N+2) }\Bigl( \frac{v}{|B|}
\Bigr)^{\frac{1}{N}}   +o(v^{\frac{1}{N}})
\ee
as $v=\left| B_g(y_0,r)\right|_g \to 0.$
\end{cor}

\begin{pf}
By the volume expansion (\ref{eq:expvolBgr}) of geodesic balls we have 
\begin{align*}
\frac{1}{r} \Bigl(\frac{v}{
  |B|}\Bigr)^{\frac{1}{N}}=\Bigl(\frac{|B_g(y_0,r)|_g}{r^{N}
  |B|}\Bigr)^{\frac{1}{N}}&= 1 - \frac{1}{6N(N+2)}S(y_0)r^2+o(r^2)\\
&= 1 -
\frac{1}{6N(N+2)}S(y_0)\Bigl(\frac{v}{|B|}\Bigr)^{\frac{2}{N}}
+o\Bigl(\frac{v}{|B|}\Bigr)^{\frac{2}{N}}
\end{align*}
as $v=\left| B_g(y_0,r)\right|_g \to 0$. Combining this with
Proposition~\ref{lem:expa_un2B}, we get the result.
\end{pf}

Next, we wish to derive an expansion of $\nu_2$ on small  geodesic ellipsoids centered at $y_0 \in \M$. For this we assume in the following that the orthonormal
basis $E_i$, $i=1,\dots,N$ of $T_{y_0}\M$ is chosen such that
\begin{equation}
  \label{eq:37}
R_{ij}= 0  \qquad \text{for $i \not=j$}.
\end{equation}
Moreover we put
\begin{equation}
  \label{eq:30}
b_i=b^i:= \frac{1}{3(N+2)} (R_{ii} -\frac{S(y_0)}{N}) \qquad
\text{for $i=1,\dots,N$,}
\end{equation}
and we note that  $\sum \limits_{i=1}^N b_i= 0$ since $S(y_0)= \sum \limits_{i=1}^N
R_{ii}$. For $r>0$ small, we then consider the geodesic ellipsoids
$E(y_0,r):= F_r(B) \subset \M$, where
$$
F_r: B \to \M, \qquad F_r(x)= \textrm{Exp}_{y_0}\bigl(r (1+r^2b_i)x^i E_i\bigr).
$$
The special choice of the values $b_i$ gives rise to the following
asymptotic expansion depending only on the
scalar curvature at $y_0$.

\begin{pro}
\label{sec:expansion-mu_2-small-4}
As $r \to 0$, we have
\begin{equation}
  \label{eq:29}
\nu_2(E(y_0,r),g) =\frac{1}{r}+\frac{2r}{3N(N+2) }S(y_0)+o(r)
\end{equation}
and
\begin{equation}
  \label{eq:31}
|E(y_0,r)|_g=|B_g(y_0,r)|_g+O(r^{N+4})= r^N |B| \Bigl(1 -
\frac{1}{6(N+2)}r^2 S(y_0) + O(r^4)\Bigr).
\end{equation}
\end{pro}

\begin{pf}
We consider the pull back metric $h_r$ on $B$
of $g$ under the map $F_r$ rescaled with the
factor $\frac{1}{r^2}.$ Then we have
\begin{align}
\label{eq:exp-hr}
[h_r]_{ij}(x)&=(1+r^2 b_i)(1+r^2 b_j)[g_r]_{ij}((1+r^2 b_k)x^k
e_k)=[g_r]_{ij}(x)+ r^2(b_i+b_j)\d_{ij}+O(r^4)\\
&= \delta_{ij} + r^2 \Bigl(\frac{1}{3}R_{kilj}x^kx^l
 + 2b_i\d_{ij}\Bigr)+O(r^3)\nonumber
\end{align}
uniformly in $x \in B$ (where $g_r$ is defined as in the proof of Proposition~\ref{lem:expa_un2B}).  Setting $|h_r|= \det ([h_r]_{ij})_{ij}$, we deduce the expansion
\begin{equation}
  \label{eq:7}
|h_r|(x)= |g_r|(x) + 2 r^2 \sum_{i=1}^N b_i + O(r^4)= |g_r|(x)+ O(r^4)
\qquad \text{for $x \in B$.}
\end{equation}
This implies that
$$
|E(y_0,r)|_g = r^N |B|_{h_r}= r^N \Bigl(|B|_{g_r}+O(r^4)\Bigr)
=|B_g(y_0,r)|_g+O(r^{N+4}),
$$
as claimed in (\ref{eq:31}).  We now turn  to \eqref{eq:29}. For this we first note that, denoting by $(h_r^{ij})_{ij}$ the inverse of the matrix $([h_r]_{ij})_{ij}$, we have
\begin{equation}
  \label{h-inverse-exp}
 h_r^{ij}(x)= \delta^{ij} - r^2 \Bigl(\frac{1}{3}R_{k}\!{}^i\!{}_l{}^j x^kx^l
 + 2 b_i \d_{ij}\Bigr)+O(r^3)
\end{equation}
by (\ref{eq:exp-hr}), whereas (\ref{eq:3}) and (\ref{eq:7}) yield
\begin{equation}
  \label{eq:23}
 \frac{\de}{\de x^i}\sqrt{|h_r|}= -\frac{r^2}{3}R_{ki}x^k+O(r^3) \qquad \text{for $i=1,\dots,N$.}
\end{equation}
Moreover, since $\nu_2(B,h_r)=r
\nu_2(E(y_0,r),g)$, the asserted expansion (\ref{eq:29}) is equivalent to
\begin{equation}
  \label{eq:33}
\nu_2(B,h_r) =1+  \frac{2r^2}{3N(N+2)}S(y_0)+o(r^2).
\end{equation}
Let $\Phi_r$ be an eigenfunction for $\nu_2(B, h_r)$, normalized
such that $\int_{\de B}\Phi_r^2\,dv_{h_r}=1$ with $dv_{h_r}=
\sqrt{|h_r|}dx$. Then we have
 $$
\D_{h_r} \Phi_r=0\quad \text{in $B$},\qquad \la \n \Phi_r,\eta
\ra_{h_r}=\nu_2(B,h_r ) \,\Phi_r\quad \text{on $\de B$,}
$$
where
$$
\D_{h_r}\Phi_r=  \frac{1}{{\sqrt{|h_r|}}}\frac{\de }{\de x^i}\left(\sqrt{|h_r|}h_r^{ij} \frac{\de \Phi_r}{\de x^j}\right).
$$
Since $h_r$ converges to the Euclidean metric in $B$, the
variational characterization of $\nu_2$ implies that $\nu_2(B,h_r
)\to \nu_2(B)=1$. Moreover, as in the proof of
Proposition~\ref{lem:expa_un2B} we have $\Phi_{r_k}\to \Phi$ in
$H^1(B)$ along a sequence $r_k \to 0$ with some function $\Phi \in
C^2_{loc}(B)\cap C^1(\ov{B})$ satisfying
$$
\D \Phi=0\quad \text{in $B$} ,\qquad  \la \n \Phi,\eta
  \ra= \Phi\quad \text{on $\de B$},\quad
  \text{and}\quad \int_{\de B} \Phi^2\,d\sigma=1.
$$
Hence there exists a vector $a=(a_1,\dots,a_{N})=(a^1,\dots,a^{N}) \in \R^{N}$
with $|a|=1$ and such that
$$
\Phi(x)=\frac{a \cdot x}{\sqrt{|B|}}   \qquad \text{for $x
  \in \overline B$.}
$$
Again, for matters of convenience, we write $r$ instead of
$r_k$ in the following. By multiple integration by parts, using
(\ref{h-inverse-exp}) and (\ref{eq:23}), we have
\begin{align}
\nu_2(B, h_r)\, &\int_{\de B} \Phi \Phi_r d\sigma_{h_r}=
 \int_{B}\sqrt{|h_r|}h_r^{ij} \frac{\de \Phi}{\de x^i}
                                          \frac{\de \Phi_r}{\de x^j}
                                          \,dx \nonumber \\
=& \int_{B} \nabla \Phi_r \nabla \Phi  dv_{h_r}  -r^2 \int_{B}\Bigl(\frac{1}{3} R_k\!{}^i\!{}_l{}^jx^kx^l\frac{\de \Phi_r}{\de x^i} \frac{\de \Phi_r}{\de x^j}  + 2 b^i \frac{\de \Phi}{\de x^i}
                                           \frac{\de \Phi_r}{\de
                                             x^i}\Bigr)\,dx+O(r^3)\nonumber
                                           \\
=&\int_{\de B}\Phi\Phi_r\,d\sigma_{h_r} - \int_{B}\Phi_r \nabla \sqrt{|h_r|}
                                  \cdot \nabla \Phi \,dx
 \nonumber\\
&-r^2 \int_{B}\Bigl(\frac{1}{3} R_k\!{}^i\!{}_l{}^jx^kx^l\frac{\de \Phi_r}{\de x^i} \frac{\de \Phi_r}{\de x^j}  + 2 b^i \frac{\de \Phi}{\de x^i}
                                           \frac{\de \Phi_r}{\de
                                             x^i}\Bigr)\,dx
                                          +O(r^3)
                                          \nonumber \\
=&\: \int_{\de
B}\Phi\Phi_r\,d\sigma_{h_r}\nonumber\\
&+r^2 \int_{B}\Bigl(\frac{\Phi_r}{3}  R_{i}{}^j x^i \frac{\de \Phi}{\de x^j}-
\frac{1}{3} R_k\!{}^i\!{}_l{}^jx^kx^l\frac{\de \Phi_r}{\de x^i} \frac{\de \Phi_r}{\de x^j}  -
 2 b^i \frac{\de \Phi}{\de x^i}
                                           \frac{\de \Phi_r}{\de
                                             x^i}\Bigr)\,dx
                                          +O(r^3)\nonumber
\end{align}
Since $ \int_{\partial B}\Phi\Phi_r\,d\sigma_{h_r}\to  1$ and $\Phi_r\to \Phi$ in
$H^1(B)$ as $r\to0$, we infer, using (\ref{eq:25-new}) and (\ref{eq:27-new}), that
\begin{align}
\nu_2(B,h_r)=&\:1+r^2  \int_{B}\Bigl(\frac{\Phi}{3}  R_{i}{}^j x^i \frac{\de \Phi}{\de x^j}-
\frac{1}{3} R_k\!{}^i\!{}_l{}^jx^kx^l\frac{\de \Phi}{\de x^i} \frac{\de \Phi}{\de x^j}  -
 2 b^i \frac{\de \Phi}{\de x^i}
                                           \frac{\de \Phi}{\de
                                             x^i}\Bigr)\,dx
                                          +o(r^2)\nonumber \\
=&\:1+r^2  \int_{B} \Bigl(   \frac{1}{3|B|}  a_k x^k R_{i}^j a_j-
\frac{1}{3|B|} R_k\!{}^i\!{}_l{}^jx^kx^l a_i a_j   -
 \frac{2}{|B|} b^i  a_i^2\Bigr) \,dx
                                          +o(r^2)\nonumber \\
=&\: {1}+{2r^2}\Bigl(\frac{1}{3(N+2)}Ric_{y_0}(A,A)
- b^i  a_i^2\Bigr) +o(r^2)\nonumber
\end{align}
with $A := a^i E_i \in T_{y_0}\M$. Combining this with (\ref{eq:37}) and (\ref{eq:30}), we conclude that
\begin{align}
\nu_2(B,h_r)=&\: {1}+{2r^2}(a^{i})^2 \Bigl(\frac{1}{3(N+2)} R_{ii} -b_i \Bigr) +o(r^2)\nonumber\\
=&\: {1}+ \frac{2r^2}{3N(N+2)} S(y_0)  +o(r^2). \nonumber
\end{align}
 Hence we have shown
(\ref{eq:33}), as required.
\end{pf}

\begin{cor}
\label{sec:expansion-mu_2-small-2} We have
\be\label{eq:expmu_2-ellipsoid}
\nu_2(E(y_0,r),g)=\Bigl(\frac{v}{|B|}\Bigr)^{-\frac{1}{N}}+\frac{S(y_0)}{2N(N+2) } \left( \frac{v}{|B|}
\right)^{\frac{1}{N}}   +o(v^{\frac{1}{N}}) \ee
as $v=\left| E(y_0,r)\right|_g \to 0.$
\end{cor}

\begin{pf}
This follows readily by combining   (\ref{eq:29}) and (\ref{eq:31}).
\end{pf}

%
\section{A local upper bound for $\nu_2$}\label{s:Ubm2}
The aim of this section is to complete the proof of Theorem~\ref{pro:Main-rslt}. We note that Theorem~\ref{pro:Main-rslt}(i) follows immediately from
Corollary \ref{sec:expansion-mu_2-small-2}, and the lower bound in
(ii) is a direct consequence of (i). Hence it remains to establish the upper bound (ii).
For this we fix $r_0>0$ less than the injectivity radius of $\M$ at
$y_0$.  Throughout this section, we consider a sequence of numbers
$r_k \in (0,\frac{r_0}{4})$ such that $r_k \to 0$ as $k \to \infty$,
and we suppose that we are given regular domains $\O_{r_k}\subset
B_g(y_0,r_k)$, $k \in \N$.  In this setting, we will show the following asymptotic upper bound.

\begin{thm}
 \label{sec:local-upper-bound-2}
We have
\begin{equation}
  \label{eq:Omega_r-exp}
 \nu_2(\Omega_{r_k},g) \le \Bigl(\frac{|\Omega_{r_k}|_g}{|B|}\Bigr)^{-\frac{1}{N}}  + \frac{S(y_0)}{2N(N-2)}
\Bigl(\frac{|\Omega_{r_k}|_g}{|B|}\Bigr)^{\frac{1}{N}} + o(|\Omega_{r_k}|_g^{\frac{1}{N}}) \qquad \text{as $k \to \infty$.}
\end{equation}
\end{thm}

This result obviously implies the upper bound in Theorem~\ref{pro:Main-rslt}(ii), so the proof of Theorem~\ref{pro:Main-rslt} is finished once we have established Theorem~\ref{sec:local-upper-bound-2}.

The remainder of this section is devoted to the proof of Theorem~\ref{sec:local-upper-bound-2}. In order to keep the notation as simple
as possible, we will write $r$ instead of $r_k$ in the following. As in the previous sections, we rescale the problem, but we first need to identify suitable center points for the rescaling procedure. For this, we need the following observation.

\begin{lem}
\label{sec:local-upper-bound-1}
There exists a point $p_r\in
B_g(y_0,2r)$ with
\be\label{eq:centmass}
 \int_{\de \O_r}
\textrm{Exp}_{p_r}^{-1}(q)\,d\sigma_g(q)=0 \; \in T_{p_r}\M.
\ee
\end{lem}

\begin{pf}
Consider the function
$$
J:  \overline{B_g(y_0,2r)} \to \R,\qquad  J(p)=  \int_{\de \O_r}| \textrm{Exp}_{p}^{-1}(q)|^2_g\,d\sigma_g(q)=
\int_{\de \O_r}\textrm{dist}_g(p,q)^2\,d\sigma_g(q).
$$
Since $r<r_0$ and $\O_r \subset B_g(y_0,r)$, the function $J$ is differentiable with
$$
dJ(p)[v]= -2 \int_{\de \O_r}\la
\textrm{Exp}_{p}^{-1}(q),v\ra_g\,d\sigma_g(q)\qquad \text{for all $v\in T_p\M.$}
$$
Since $J(y_0) \le r^2 \sigma_g(\partial \Omega_r)$ and
$$
J(p) \ge r^2\sigma_g(\partial \Omega_r)\quad \text{for $p \in \partial B_g(y_0,2r)$},
$$
there  exists a point $p_r\in
B_g(y_0,2r)$ with  $J(p_r) = \min\{J(p): p \in B_g(y_0,2r)\}$. Hence $p_r$ is a critical point of $J$,  and this implies \eqref{eq:centmass}.
\end{pf}

Next we note that,  for $r>0$ small enough, we have $|B_g(p_r, 2 r)|_g >|B_g(y_0,r)|_g$, and thus there exists a unique $\rho_r \in (0,2r)$ with
$$
| \O_r|_g=| B_g(p_r,\rho_r) |_g.
$$
Since $p_r \to y_0$ as $r \to 0$, we have, similarly as in \eqref{eq:expvolBgr}, the volume expansion
$$
 \frac{| \O_r|_g}{|B|}=\frac{| B_g(p_r,\rho_r) |_g}{|B|} = \rho_r^N  \Bigl(1-\frac{S(y_0)}{6(N+2)}\rho_r^2 +o(\rho_r^2)\Bigr)
$$
and thus
\begin{equation}
  \label{volume-new}
\Bigl( \frac{| \O_r|_g}{|B|}\Bigr)^{\frac{1}{N}}= \rho_r
 \Bigl(1-\frac{S(y_0)}{6N(N+2)}\rho_r^2 +o(\rho_r^2)\Bigr).
\end{equation}
Consequently, Theorem~\ref{sec:local-upper-bound-2} is proved once we establish the following:
\begin{equation}
  \label{eq:25-new-0}
\nu_2(\Omega_r,g) \le \frac{1}{\rho_r} + \frac{2\rho_r}{3N(N+2)}S(y_0) + o(\rho_r) \qquad \text{as $r \to 0$.}
\end{equation}
We now consider a rescaled version of (\ref{eq:25-new-0}). For this we note that
$$
B_g(p_r,\rho_r) \subset B_g(p_r,2r) \subset B_g(y_0,4r) \qquad \text{and}\qquad \Omega_r \subset B_g(y_0,r) \subset B_g(p_r,3r),
$$
and we let
$$
y \mapsto E_i^y \in T_y \M,\qquad i=1,\dots,N
$$
denote a smooth orthonormal frame on $B_g(y_0,r_0)$. We consider the maps
$$
\Psi_r: \R^{N} \to \M,\qquad \Psi_r(x)=
\Psi(\rho_rx)=\textrm{Exp}_{p_r}(\rho_r x^i E_i^{p_r}).
$$
Moreover, we set
 \be\label{eq:defUr}
B^r:= \frac{3r}{\rho_r}B \qquad \text{and}\qquad
U_r:=\Psi_r^{-1} (\O_r) \subset B^r, \ee
and we consider the pull back metric of $g$ under the map $B^r \to
\M,\; x \mapsto \Psi_r(\rho_r x)$, rescaled with the factor
$\frac{1}{\rho_r^2}$. We denote this metric on $B^r$ by $g_r$, and
we point out that this definition differs from the notation used in
the proof of Proposition~\ref{lem:expa_un2B}. Nevertheless, since
$\dist(p_r,y_0)=O(r)$, we have, in $C^1$-sense,
\begin{equation*}
\la R_y(E_i^{p_r},E_j^{p_r})E_k^{p_r},E_l^{p_r} \ra =  R_{ijkl}+O(r)  \quad \text{as $r \to 0\quad$ with}\quad
R_{ijkl}:= \la R_{y_0}(E_i^{y_0},E_j^{y_0})E_k^{y_0},E_l^{y_0} \ra
\end{equation*}
for $i,j,k,l=1,\dots,N$. We also set $R_{ij}:= Ric_{y_0}(E_i^{y_0},E_j^{y_0})$. As in Section~\ref{s:pn}, we freely vary the (upper or lower) position of the indices of $R_{ijkl}$ and $R_{ij}$ without changing the value of these constants.  We then infer from (\ref{exp-fund-1}) and (\ref{exp-fund-2}) that
\be\label{eq:metexp}
\begin{array}{rllll}
\displaystyle (g_r)_{ij}(x)&=\delta_{ij}+\frac{\rho_r^2}{3}
R_{kilj} x^kx^l + O(r \rho_r^2) ,\\[3mm]
\displaystyle g_r^{ij}(x)&=\delta^{ij}-\frac{\rho_r^2}{3}
R_{k}\!{}^i\!{}_l{}^j x^kx^l + O(r \rho_r^2) ,\\[3mm]
 \displaystyle dv_{g_r}(x)=\sqrt{|g_r(x)|}\,dx&=\left(
   1-\frac{\rho_r^2}{6}\, R_{lk} x^lx^k +O(r \rho_r^2) \right)dx,\\[3mm]
\frac{\partial}{\partial x_i} \sqrt{|g_r(x)|} &= - \frac{\rho_r^2}{3} R_{ki}x^k+ O(r\rho_r^2)
\end{array}
\ee
uniformly on $B^r$ as $r \to 0$, where $(g_r^{ij})_{ij}$ denotes the inverse of the matrix
$([g_r]_{ij})_{ij}$ and $|g_r|$ is the
determinant of $g_r$.  In particular
\be
\label{eq:metexp-1-new}
(g_r)_{ij}(x)=\delta_{ij}+O(r^2) \quad \text{and}\quad
dv_{g_r}(x)=(1+O(r^2))dx \qquad \text{uniformly on $B^r$.}
\ee
Moreover, by the same arguments as in the proof of Lemma~\ref{sec:local-expans-nu_2-surface} we have
\begin{equation}
  \label{eq:boundary-integral-modified}
\int_{\de U_r}f(x)  \,d\sigma_{g_{r}}= (1+O(\rho_r^2)) \int_{\de U_r}f(x) d\s,
\end{equation}
for every $f\in C^1(\de U_r)$. Moreover,
\begin{equation}
   \label{exp|g_r|-surface-modified}
\int_{\de B}f(x)  \, d\sigma_{g_r}(x)= \int_{\de B} \Bigl(1-\frac{\rho_r^2}{6}R_{lk}x^l x^k\Bigr)f(x)d\s  + O(r \rho_r^2 )\int_{\de B}f(x) d\s
\end{equation}
for every $f \in C^1(\de B)$. Here, similarly as in Lemma~\ref{sec:local-expans-nu_2-surface}, the bounds for the terms $O(\rho_r^2)$ and $O(r\rho_r^2)$ are uniform in $f$. Observe also that $\nu_2(U_r,g_r)=\frac{\nu_2(\O_r,\,g)}{\rho_r}$, so that
(\ref{eq:25-new-0}) is equivalent to
\begin{equation}
  \label{eq:25-new_1}
\nu_2(U_r,g_r) \le 1 + \frac{2\rho_r^2}{3N(N+2)}S(y_0) + o(\rho_r^2) \qquad \text{as $r \to 0$.}
\end{equation}
The remainder of this section will be devoted to the proof of (\ref{eq:25-new_1}). By construction we have $|U_r|_{g_r}={\rho_r}^{-N}|\Omega_r|_{g}={\rho_r}^{-N}|B_g(p_r,\rho_r)
|_g= |B|_{g_r}$, and thus
 \be\label{eq:volUdx}
|U_r|_{g_r}=|B|_{g_r}=(1+O(r^2))|U_r|=(1+O(r^2))|B|
\ee
by \eqref{eq:metexp-1-new} and the fact that $U_r \subset B^r$ and $B \subset B^r$. Setting
$$
{ f}_i: \R^{N} \to \R,\qquad f_i(x)={x^i},
$$
we also find that $\int_{\de U_r}f_id\sigma_{g_r}=0$ for $i=1,\dots,N$ by
\eqref{eq:centmass}. Moreover,
$$
\int_{U_r}|\n f_i|_{g_r}^2dv_{g_r} = \int_{U_r}  g_r^{jk} \frac{\partial f_i}{\partial x^j}
\frac{\partial f_i}{\partial x^k} dv_{g_r} = \int_{U_r} g_r^{ii} d v_{g_r}\qquad \text{for $i=1,\dots,N$.}
$$
Hence the variational characterization of
$\nu_2$ yields \be\label{eq:mudlesSm}
  \nu_2(U_r,g_r)\leq\frac{\displaystyle
\sum_{i=1}^N\int_{U_r}|\n f_i|_{g_r}^2dv_{g_r}  }{\displaystyle
\sum_{i=1}^N \int_{\de U_r}f_i^2d\sigma_{g_r}}=  \frac{\displaystyle
\sum_{i=1}^N\int_{U_r} g_r^{ii}   dv_{g_r}
}{\displaystyle \int_{\de U_r}|x|^2d\sigma_{g_r}} .
 \ee
In the following, $B_2:= 2B \subset \R^N$ denotes the euclidean ball centered at the origin with radius 2. Moreover, we let $|U_r\triangle B| = |U_r \setminus B|+ |B \setminus U_r|$ denote the symmetric distance of the sets $U_r$ and $B$ with respect to the standard Lebesgue measure on $\R^N$.

\begin{lem}\label{lem:appest}
In the above setting, we have
\be\label{eq:estintUrgsq}
 \int_{\de U_r}|x|^2\,d\sigma_{g_r}\geq N|B|-\frac{|B|\rho^2_r}{6}S(y_0)+O(\rho_r^2|U_r\triangle B|)+\frac{N+1}{4} |U_r\setminus B_2
 | _{g_r} +o(\rho_r^2),
 \ee
and
 \be\label{eq:estintUrgsq1}
\sum_{i=1}^N\int_{U_r} g_r^{ii}   dv_{g_r}   =N|B|
-\frac{(N-2)|B|}{6(N+2)}\rho_r^2S(y_0)+O(\rho_r^2|U_r\triangle B|)+
O(r^2|U_r\setminus B_2 |_{g_r})+o(\rho_r^2)
  \ee
as $r \to 0$.
\end{lem}
\begin{pf}
We first note that, by \eqref{eq:metexp-1-new}, the symmetric distance $|U_r\triangle B|_{g_r}: = |U_r \setminus B|_{g_r}+ |B \setminus U_r|_{g_r}$ with respect to the metric $g_r$ satisfies
\begin{equation}
  \label{exp-sym-dist}
|U_r\triangle B|_{g_r} = (1+O(r^2))|U_r\triangle B|.
\end{equation}
Next we consider the $C^1$-vector field $V: B^r \to \R, \: V(x)= |x|x$. Using \eqref{eq:metexp}, we have
\begin{align}
 G(x)&:= \textrm{div}_{g_r}(V)=\frac{1}{ \sqrt{|g_r|}(x)}
\frac{\partial}{\partial x^i}\left(|x|x^i \sqrt{|g_r|}(x)  \right)=|x| \Bigl[(N+1) -  \frac{x \cdot \nabla \sqrt{|g_r|(x)}}{\sqrt{|g_r|(x)}} \Bigr]\nonumber\\
 \nonumber&=|x| \Bigl[(N+1) -  \frac{\frac{\rho_r^2}{3}Ric_{y_0}(X,X)+ O(r \rho_r^2)}{\sqrt{|g_r|(x)}} \Bigr]
= {|x|} \Bigl[(N+1)  -\frac{\rho_r^2}{3}Ric_{y_0}(X,X)+
O(r \rho_r^2) \Bigr]\\
&={|x|} \Bigl[(N+1)  -\frac{\rho_r^2}{3}Ric_{y_0}(X,X)+ O(r\rho_r^2 ) \Bigr]\qquad \text{uniformly for $x \in B^r$ as $r \to 0$.}\label{eq:divVr}
\end{align}
In particular, $G(x)=|x|\Bigl[(N+1) + O(r^2)\Bigr]$ for $x \in B^r$ as $r \to 0$, so for $r>0$ sufficiently small we have
\begin{equation}
  \label{eq:9-1}
  G(x) \ge \frac{N+1}{2} |x|  \qquad \text{for $x \in B^r$.}
\end{equation}
We also recall that, as a consequence of Gauss' Lemma (see e.g. \cite[Corollary 5.2.3]{jost:RiemannianGeometryandGeometricAnalysis}), the unit outer normal on $\partial B$ with respect to the metric $g_r$ is simply given by $\eta_r(x)=x$ for every small $r>0$ and $x \in \partial B$. Using the divergence formula with respect to the metric $g_r$ and (\ref{exp|g_r|-surface-modified}), we therefore find that
\begin{equation}
  \int_{B}  G dv_{g_r}\!=\!\int_{\de B}  d\sigma_{g_r}\!=\!
\int_{\de B}  \Bigl(1 - \frac{\rho_r^2}{6} R_{lk}x^lx^k \Bigr)\,d\sigma   +O(r \rho_r^2)
= N|B| -\frac{|B| \rho_r^2}{6} S(y_0) +o(\rho_r^2).\label{eq:intB}
\end{equation}
Here we used the fact that $\int_{\de B} x^lx^k\,d\sigma = \delta^{lk} \frac{\sigma(\partial B)}{N} = \delta^{lk}|B|$ in the last step. Moreover, using again that $|x|_{g_r} = |x|$ and thus $|V(x)|_{g_r}=|x|^2$ for $x \in B^r$ by Gauss' Lemma, we find that
 \be\label{eq:intUr}
  \int_{U_r} G\,
dv_{g_r}=\int_{\de U_r}\la V ,\eta_r \ra_{g_r} d\sigma_{g_r}\leq
\int_{\de U_r}|V|_{g_r} \,d\sigma_{g_r}= \int_{\de U_r}|x|^2 \,d\sigma_{g_r},
\ee
where $\eta_r $ is the outer unit normal of $\de U_r$ with respect to $g_r$. Next we estimate
\begin{align}
\int_{U_r}G\,dv_{g_r}&- \int_{B} G\,dv_{g_r}=
\int_{U_r\setminus B} G\,dv_{g_r}-\int_{B\setminus U_r}G\,dv_{g_r}\nonumber\\
&=\int_{U_r\setminus B}(1-\frac{1}{|x|})G\,dv_{g_r}+\int_{U_r\setminus
B}\frac{G}{|x|}\,dv_{g_r}  -\int_{B \setminus U_r} G\,dv_{g_r} \nonumber\\
&\geq  \int_{U_r \setminus  B}(1-\frac{1}{|x|})G\,dv_{g_r} +\int_{U_r\setminus B}\frac{G}{|x|}\,dv_{g_r} -\int_{B \setminus U_r} \frac{G}{|x|}\,dv_{g_r}.\label{eq:12}
\end{align}
Here we note that, by (\ref{eq:9-1}),
\begin{equation}
  \label{eq:6}
  \int_{U_r \setminus  B}(1-\frac{1}{|x|})G\,dv_{g_r} \ge   \int_{U_r \setminus  B_2}(1-\frac{1}{|x|})G\,dv_{g_r}
 \ge   \int_{U_r \setminus  B_2}\frac{G}{|x|}\,dv_{g_r} \ge   \frac{N+1}{2} |U_r \setminus B_2|_{g_r}
\end{equation}
 and, by \eqref{eq:volUdx} and (\ref{eq:divVr}),
\begin{align}
\int_{U_r\setminus B}\frac{G}{|x|}\,dv_{g_r} &-\int_{B \setminus U_r} \frac{G}{|x|}\,dv_{g_r}=
\int_{U_r\setminus B}\Bigl(\frac{G}{|x|}-(N+1)\Bigr)\,dv_{g_r} -\int_{B \setminus U_r}\Bigl(\frac{G}{|x|}-(N+1)\Bigr)\,dv_{g_r}\nonumber\\
&=  \frac{\rho_r^2}{3} \Bigl[\int_{B \setminus U_r} \Bigl(R_{ij}x^ix^j+ O(r) \Bigr)\,dv_{g_r}
- \int_{U_r\setminus B}\Bigl(R_{ij}x^ix^j+ O(r) \Bigr)\,dv_{g_r}\Bigr] \nonumber\\
&=  \frac{\rho_r^2}{3} \Bigl[\int_{B \setminus U_r} R_{ij}x^ix^j\,dv_{g_r}
- \int_{U_r\setminus B} R_{ij}x^ix^j\,dv_{g_r} + O(r |B \triangle U_r|_{g_r}) \Bigr] \nonumber\\
&=  \frac{\rho_r^2}{3} \Bigl[- \int_{U_r\setminus B_2}R_{ij}x^ix^j \,dv_{g_r} + O(|B \triangle U_r|_{g_r}) \Bigr] \nonumber\\
&=  O( r^2 |U_r \setminus B_2
|_{g_r})  + O(\rho_r^2 |U_r \triangle B|) +o(\rho_r^2) \label{eq:13},
\end{align}
where in the last step we used (\ref{exp-sym-dist}) and the fact that
$U_r \subset B^r$. Combining (\ref{eq:12}),~(\ref{eq:6}) and (\ref{eq:13}),  we obtain that
$$
\int_{U_r}G\,dv_{g_r}\geq \int_{B} G\,dv_{g_r}
+O(\rho_r^2|U_r\triangle B|) + \frac{N+1}{4}|U_r \setminus B_2|_{g_r} \qquad \text{for $r>0$ sufficiently small.}
$$
 Combining this with \eqref{eq:intB} and
\eqref{eq:intUr}, we get the inequality
\be
\int_{\de U_r}|x|^2 \,d\sigma_{g_r}\geq  N|B| - \frac{|B|\rho_r^2}{6} S(y_0) +\frac{N+1}{4} |U_r\setminus B_2 |_{g_r}
+O(\rho_r^2|U_r\triangle B|)+o(\rho_r^2)
\ee
as $r \to 0$, which is \eqref{eq:estintUrgsq}. Next, using \eqref{eq:volUdx}, we estimate similarly as in (\ref{eq:13}),
\begin{align*}
&\sum_{i=1}^N\int_{U_r}g_r^{ii}
dv_{g_r}=\int_{U_r}[N+\frac{\rho_r^2}{3}R_{ij}x^i x^j+ O(r \rho_r^2)]dv_{g_r}\\
&=
\int_{B}(N+\frac{\rho_r^2}{3}R_{ij}x^i x^j+o(\rho_r^2))dv_{g_r}+  \frac{\rho_r^2}{3}
\int_{U_r \setminus B_2} \Bigl(R_{ij}x^i x^j +O(r)\Bigr)dv_{g_r} +O(\rho_r^2|U_r\triangle
B|_{g_r})\\
&=
N|B|  -\frac{(N-2)|B|}{6(N+2)}\rho_r^2S(y_0) +o(\rho_r^2)+O(r^2 |U_r \setminus B_2|_{g_r}) +O(\rho_r^2|U_r\triangle
B|_{g_r})\\
&=
N|B|  -\frac{(N-2)|B|}{6(N+2)}\rho_r^2S(y_0) +O(r |U_r \setminus B_2|_{g_r}) +O(\rho_r^2|U_r\triangle
B|)+o(\rho_r^2)
\end{align*}
as $r \to 0$, as claimed in \eqref{eq:estintUrgsq1}.

\end{pf}\\

We may now complete the\\

\noindent {\em Proof of Theorem~\ref{sec:local-upper-bound-2}.} As noted before, it suffices to prove (\ref{eq:25-new_1}), since (\ref{eq:25-new_1}) is equivalent to (\ref{eq:25-new-0}) and (\ref{eq:25-new-0}) is equivalent to (\ref{eq:Omega_r-exp}) by the volume expansion (\ref{volume-new}). To prove (\ref{eq:25-new_1}) for $r=r_k \to 0$ as $k \to \infty$, we argue by contradiction and assume that there exists $\e_0>0$  and a subsequence -- still denoted by $(r_k)_k$ -- such
 \begin{equation}
   \label{contradiction-25-new-1}
\nu_2(U_{r_k},g_{r_k}) \ge 1 + \Bigl(\frac{2}{3N(N+2)}S(y_0) +\e_0\Bigr)\rho_{r_k}^2  \qquad \text{for all $k \in \N$.}
\end{equation}
We first claim that
\begin{equation}
  \label{eq:22}
 |U_{r_k} \setminus B_2|_{g_{r_k}}  =O(\rho_{r_k}^2) \qquad \text{as $k \to \infty$.}
\end{equation}
Indeed, if, by contradiction, for a subsequence we have $ |U_{r_k} \setminus B_2|_{g_{r_k}} \ge k\rho_{r_k}^2$, then the expansions \eqref{eq:estintUrgsq} and \eqref{eq:estintUrgsq1} yield that
$$
\int_{\de U_{r_k}}|x|^2d\sigma_{g_{r_k}} \ge \sum_{i=1}^N \int_{U_{r_k}} g_{r_k}^{ii}   dv_{g_{r_k}} \qquad \text{for $k$ sufficiently large.}
$$
Here we also used the fact $|U_{r_k} \triangle B|$ remains bounded as a consequence of \eqref{eq:volUdx}. Now \eqref{eq:mudlesSm} implies that $\nu_2(U_{r_k},g_{r_k}) \le  1$ for $k$ sufficiently large, contrary to (\ref{contradiction-25-new-1}).  Hence (\ref{eq:22}) is true. From (\ref{eq:22}) and \eqref{eq:estintUrgsq1} it then follows that
 \be\label{eq:estintUrgsq2}
\sum_{i=1}^N \int_{U_{r_k}} g_{r_k}^{ii}   dv_{g_{r_k}}   =N|B| + O({r_k}^2).
  \ee
 Since also $\nu_2(U_{r_k},g_{r_k}) \ge 1 +O({r_k}^2)$ by (\ref{contradiction-25-new-1}), it follows from \eqref{eq:mudlesSm} that
$$
\int_{\de U_{r_k}}|x|^2 d\sigma_{g_{r_k}}  \le N|B| +O({r_k}^2) .
$$
and thus also
\begin{equation}
  \label{eq:18}
\int_{\de U_{r_k}}|x|^2d\sigma  \le N|B| +O({r_k}^2)
\end{equation}
as a consequence of (\ref{eq:boundary-integral-modified}). On the other hand, \cite[Theorem B]{BDR} implies that
\begin{equation}
  \label{eq:19}
 \int_{\de U_{r_k}}|x|^2d\sigma \ge N|B| + \beta \Bigl(\frac{|U_{r_k} \triangle B|}{|U_{r_k}|} \Bigr)^2
\end{equation}
with a positive constant $\beta>0$. Hence, by \eqref{eq:volUdx} and (\ref{eq:18}),
\begin{equation}
 \label{est-symm-dist}
|U_{r_k} \triangle B| = O(|U_{r_k}| {r_k}^2) =O({r_k}^2) \qquad \text{as $k  \to \infty$.}
\end{equation}
Inserting this in \eqref{eq:estintUrgsq} gives
\begin{equation}
  \label{eq:final-1}
\int_{\de U_{r_k}}|x|^2d\sigma_{g_{r_k}}  \ge N|B| - \frac{|B|}{6}S(y_0) \rho_{r_k}^2 + o(\rho_{r_k}^2).
\end{equation}
Moreover, inserting (\ref{eq:22}) and (\ref{est-symm-dist}) in \eqref{eq:estintUrgsq1} gives
\begin{equation}
  \label{eq:final-2}
\sum_{i=1}^N \int_{U_{r_k}} g_{r_k}^{ii}   dv_{g_{r_k}}   =N|B| - \frac{(N-2)|B|}{6(N+2)}S(y_0) \rho_{r_k}^2 + o(\rho_{r_k}^2).
\end{equation}
Combining (\ref{eq:final-1}),~(\ref{eq:final-2}) and \eqref{eq:mudlesSm} finally yields
$$
\nu_2(U_{r_k},g_{r_k}) \le 1 + \frac{2\rho_{r_k}^2}{3N(N+2)}S(y_0) + o(\rho_{r_k}^2) \qquad \text{as $k \to \infty$,}
$$
contrary to (\ref{contradiction-25-new-1}). The proof is finished.


\section{Precise global asymptotics in the two-dimensional case}
\label{sec:global-upper-bound}
In this section we give the proof of Theorem~\ref{sec:global-upper-bound-5-0}. We shall see that most of the argument works for $N\geq 2$  except at the end of the proof of Lemma \ref{sec:global-upper-bound-3} below where we had to assume that $N=2$.  For convenience, we repeat the statement of the theorem.

\begin{thm}
\label{sec:global-upper-bound-5}
Let $(\M,g)$ be a closed Riemannian surface.  Then we have
$$
{W\!B}_{\M}(v) \;= \;\Bigl(\frac{v}{\pi}\Bigr)^{-\frac{1}{2}} +\frac{S_{\M}}{16} \Bigl(\frac{v}{\pi}\Bigr)^{\frac{1}{2}}
+o(v^{\frac{1}{2}})\qquad \text{as $v \to 0$, }
$$
where $S_{\M}$ denotes the maximum of the scalar curvature function $S$ on $\M$.
\end{thm}

The remainder of this section is devoted to the proof of this result. In view of Theorem~\ref{pro:Main-rslt}(i) and the remarks after this theorem, we only need to prove that
\begin{equation}
  \label{eq:dim2-asympt}
{W\!B}_{\M}(v) \;\le \;\Bigl(\frac{v}{\pi}\Bigr)^{-\frac{1}{2}} +\frac{S_{\M}}{16}\Bigl(\frac{v}{\pi}\Bigr)^{\frac{1}{2}}
+o(v^{\frac{1}{2}}) \qquad \text{as $v \to 0$.}
 \end{equation}
 We argue by contradiction and assume that there exists $\bar \e >0$ and a sequence of regular domains
 $\Omega_k \subset \M$ such that $v_k := |\Omega_k| \to 0$ as $k \to \infty$ and
\begin{equation}
  \label{eq:dim2-asympt-1}
\nu_2(\Omega_k,g)\;\ge \;\Bigl(\frac{v_k}{\pi}\Bigr)^{-\frac{1}{2}} +\Bigl[\frac{S_{\M}}{16}+\bar \e\Bigr]\Bigl(\frac{v_k}{\pi}\Bigr)^{\frac{1}{2}}
\qquad \text{for every $k \in \N$.}
 \end{equation}
We will show that
\begin{equation}
  \label{eq:dim2-20}
  \diam(\Omega_k) \to 0 \qquad \text{as $k \to \infty$.}
\end{equation}
Once this fact is established, we arrive at a contradiction as follows. By the compactness of $\M$, there exists $y_0 \in \M$ such that, after passing to a subsequence,
$$
\text{for every $r >0$ there exists $k_r \in \N$ such that $\Omega_{k} \subset B_g(y_0,r)$ for $k \ge k_r$.}
$$
Fix $\e< \bar \e$, and let $r_\e$ be given by Theorem~\ref{pro:Main-rslt}(ii) corresponding to these choices of $y_0$ and $\e$. Then, for $k \ge k_{r_\e}$, we have
$$
 \nu_2(\Omega_k,g) \le   {W\!B}_{\M}(v_k) \le \Bigl(\frac{v_k}{\pi}\Bigr)^{-\frac{1}{2}} +\Bigl[\frac{S(y_0)}{16}+\e\Bigr]\Bigl(\frac{v_k}{\pi}\Bigr)^{\frac{1}{2}}\le \Bigl(\frac{\pi}{v_k}\Bigr)^{\frac{1}{2}} +\Bigl[\frac{S_\M}{16}+\e\Bigr]\Bigl(\frac{v_k}{\pi}\Bigr)^{\frac{1}{2}}
$$
as a consequence of the upper estimate in Theorem~\ref{pro:Main-rslt}(ii). This contradicts (\ref{eq:dim2-asympt-1}), since $\e< \bar \e$, and thus the proof of Theorem~\ref{sec:global-upper-bound-5} is finished.
Hence it remains to prove (\ref{eq:dim2-20}), and the remainder of this section is devoted to this task. Since $\M$ is closed, it is easy to see that there exists a number $K>0$ such that
\begin{equation}
\label{sec:global-upper-bound-2}
\text{for every $r>0$, $p \in \M$ there exist $p_1,\dots,p_K \in \M$ with $B_g(p,4r) \subset \bigcup \limits_{i=1}^{K} B_g(p_k,r)$.}
\end{equation}
To prove (\ref{eq:dim2-20}), we now argue by contradiction and assume that there exists $d>0$ such that, after passing to a subsequence, $\diam (\Omega_k) \ge d$ for all $k \in \N$. In the following, we let $r_{\M}$ denote the injectivity radius of $\M$, and we put $r_0:=\min \{\frac{r_{\M}}{5}, \frac{d}{7}\}$. We also let $\phi \in C^\infty_0(\R,\R)$  be a function such that
$$
\text{$\phi \equiv 2$ on $(-\infty,0]$} ,\; \phi(r_0^2) =1, \; \phi((2r_0)^2) = \frac{1}{K+1},\;  \text{$\phi' <0$ on $(0,(4r_0)^2)$ and $\phi \equiv 0$ on $[(4r_0)^2,\infty)$},
$$
where $K$ has the property in~\eqref{sec:global-upper-bound-2}. For $p \in \M$ and $v \in T_p \M$ we consider the function
$$
f_{p,v} \in C^\infty(\M),\qquad f_{p,v}(q)= \left \{
  \begin{aligned}
& \phi' (\dist(p,q)^2) \la \textrm{Exp}_{p}^{-1}(q), v \ra, &&\qquad q \in B_g(p, 4r_0);\\
&0, &&\qquad q \not \in B_g(p,4r_0).
 \end{aligned}
\right.
$$
Since $\M$ is compact, we find that
\begin{equation}
  \label{eq-dim-2:11}
c_0 :=  \sup \{ |\nabla f_{p,v}(q)|_g\::\: q,p \in \M, \:v \in T_p \M,\: |v|_g=1\} < \infty.
\end{equation}

\begin{lem}
\label{sec:global-upper-bound-3}
There exists points $p_k \in \M$ and vectors $v_k \in T_{p_k}\M$, $k \in \N$ with $|v_k|_g=1$ and the following properties:
\begin{itemize}
\item[(i)] $\partial \Omega_k  \cap B_g(p_k,2r_0) \not= \varnothing$ for all $k \in \N$.
\item[(ii)] Setting $f_k:= f_{p_k,v_k} \in C^\infty(\M)$, we have
$\int_{\partial \Omega_k} f_k\,d\sigma_g = 0$ for all $k \in \N$.  Moreover,
\begin{equation}
  \label{replace-iii}
 c_1:= \liminf_{k \to \infty} \int_{\partial \Omega_k} f_k^2 \,d\sigma_g>0.
 \end{equation}
 \end{itemize}
\end{lem}

\begin{pf}
We fix $k \in \N$ and consider the functional
$$
J: \M \to \R, \qquad J(p)= \int_{\partial \Omega_k} \phi(\dist(p,q)^2)\,d\sigma_g(q).
$$
Since $\M$ is compact, there exists a point $p_k \in \M$ such that $J(p_k)= \max \limits_{\M}J$. We claim that
\begin{equation}
  \label{eq:dim2-15}
\dist(p,\partial \Omega_k) < 2 r_0.
\end{equation}
Indeed, suppose by contradiction that $\dist(p_k,\partial \Omega_k) \ge 2 r_0$. Then
$$
J(p_k)  \le \phi ([2r_0]^2)\, \sigma_g\Bigl(\partial \Omega_k \cap [B_g(p_k,4 r_0) \setminus B_g(p_k,2 r_0)]\Bigr) \le  \frac{\sigma_g(\partial \Omega_k \cap B_g(p_k,4 r_0))}{K+1}.
$$
On the other hand, by (\ref{sec:global-upper-bound-2}) there exists a point $\bar p \in \M$ such that $\sigma_g(\partial \Omega_k \cap B_g(\bar p, r_0)) \ge \frac{\sigma_g(\partial \Omega_k \cap B_g(p_k,4 r_0))}{K}$,
 and thus
$$
J(\bar p) \ge \phi(r_0^2)\sigma_g(\partial \Omega_k \cap B_g(\bar p, r_0)) \ge \frac{\sigma_g(\partial \Omega_k \cap B_g(p_k,4 r_0))}{K} >J(p_k),
$$
contradiction.  Hence (\ref{eq:dim2-15}) is true, and thus (i) follows.  By the maximization property of $p_k$, we have
$$
0= dJ(p_k)[v]= -2 \int_{\de \O_k}\phi' (\dist(p_k,q)^2) \la
\textrm{Exp}_{p_k}^{-1}(q),v\ra_g\,d\sigma_g = -2 \int_{\de \Omega_k} f_{p_k,v}\,d\sigma_g  \qquad \text{for all $v\in T_{p_k}\M,$}
$$
hence the first part of (ii) follows independently of the choice of $v_k$. To prove (\ref{replace-iii}) for suitable $v_k \in T_{p_k}\M$ with $|v_k|_g=1$, , we choose orthonormal vectors $v_{k_1}, v_{k_2} \in T_{p_k}\M$ (with respect to $g$) for every $k \in \N$ .
 With $\kappa: = \inf \bigl \{[\phi']^2(r)\::\: r \in [2r_0,3r_0] \bigr\}  >0$ and
$$
\Gamma_k:= \partial \Omega_k \cap [B_g({p_k},3 r_0) \setminus \overline{B_g({p_k},2 r_0)}]\qquad \text{for $k \in \N$,}
$$
we then have
 \begin{align}
\sum_{i=1}^2  \int_{\partial \Omega_k} f_{p_k,v_{k_i}}^2 \,d\sigma_g &\ge \kappa \sum_{i=1}^2 \int_{\Gamma_k} \la \textrm{Exp}_{{p_k}}^{-1}(q), v_{k_i} \ra^2\,d\sigma_g(q) = \kappa \int_{\Gamma_k}|\textrm{Exp}_{{p_k}}^{-1}(q)|_g^2\,d\sigma_g(q)\nonumber\\
&= \kappa \int_{\Gamma_k} \dist^2({p_k},q)\,d\sigma_g(q)  \ge  (2r_0)^2 \kappa \: \sigma_g(\Gamma_k). \label{align-square-int}
\end{align}
It now remains to show that
\begin{equation}
  \label{measure-estimate-2}
\liminf_{k \to \infty}\sigma_g(\Gamma_k)  > 0.
\end{equation}
Indeed, once (\ref{measure-estimate-2}) is established, we may combine it with (\ref{align-square-int}) to see that, without loss of generality,
$$
\liminf_{k \to \infty} \int_{\partial \Omega_k} f_{p_k,v_{k}}^2 \,d\sigma_g>0 \qquad \qquad \text{with $v_k:= v_{k,1}$ for $k \in \N$.}
$$
Hence (\ref{replace-iii}) holds, and the proof is then finished.  To show (\ref{measure-estimate-2}), we put $S_{r,k}:= \partial B_g(p_k,r)$ for $k \in \N$, $r>0$. Since $\diam (\Omega_k) \ge d > 6r_0$, the domain $\Omega_k$ is not contained in $B_{g}(p_k,3r_0)$.
Hence, by (i) and since $\Omega_k$ is connected, we have
\begin{equation}
  \label{eq:intersection-1}
\Omega_k \cap S_{r,k} \not=\varnothing \qquad \text{for every $r \in (2r_0,3r_0)$.}
\end{equation}
 Next, we let
$$
T_k:= \{r \in (2r_0,3r_0) \::\: S_{r,k} \subset \Omega_k\}\qquad \text{and}\qquad R_k:= (2r_0,3r_0) \setminus T_k.
$$
We claim that, after passing to a subsequence,
\begin{equation}
  \label{eq:tktozero}
|T_k| \to 0\qquad \text{as $k \to \infty$.}
\end{equation}
To see this, we may, by the compactness of $\cM$, pass to a subsequence such that $p_k \to p_0 \in \cM$ and $p_k \in B_g(p_0,r_0)$ for all $k \in \N$.  We let
$$
y \mapsto E_i^y \in T_y \M,\qquad i=1,\dots,N
$$
denote a smooth orthonormal frame on $B_g(p_0,r_0)$, and we consider the maps
$$
\Psi_k: \R^{N} \to \M,\qquad \Psi_k(x)=\textrm{Exp}_{p_k}(x^i E_i^{p_k}) \qquad \text{for $k \in \N \cup \{0\}$.}
$$
We note that $\Psi_k$ converges locally uniformly in $C^1$-sense to $\Psi_0$ as $k  \to \infty$.  Moreover, since $r_{\cM} \ge 5 r_0$,  $\Psi_k$ maps
$3 r_0 B$ diffeomorpically onto $B_g(p_0,3r_0)$ for every $k \in \N \cup \{0\}$, and there exists a constant $\alpha>1$ such that
\begin{equation}
  \label{eq:alpha-estimate-2}
\frac{1}{\alpha} \le \sqrt{\det (g_{ij}^k(x))_{ij}} \le \alpha \qquad \text{for every $x \in 3r_0  B$, $k \in \N \cup \{0\}$.}
\end{equation}
Here $g_{ij}^k$ denote the metric coefficients associated with local parametrizations $\Psi_k$, i.e.,
$$
g_{ij}^k(x)=\la d \Psi_k(x)e_i, d \Psi_k(x)e_j \ra_g \qquad \text{for $x \in \R^{N}$, $i,j =1,\dots,N$,}
$$
where $e_i \in \R^N$, $i=1,\dots,n$ denote the coordinate vectors. We set $U_k:= \Psi_k^{-1}(\Omega_k \cap  B_g({p_k},3 r_0)) $ for $k \in \N$.
Since $|\Omega_k|_g \to 0$ as $k \to \infty$, we also have that $|U_k| \to 0$ as $k \to \infty$ as a consequence of (\ref{eq:alpha-estimate-2}). Moreover, $T_k$ is given as the set of $r \in  (2r_0,3r_0)$ such that $x \in U_k$ for every $x \in \R^N$ with $|x|=r$. Hence we estimate that
$$
|U_k| \ge \omega_{N-1} \int_{2r_0}^{3r_0} r^{N-1} 1_{\text{\tiny$T_k$}}(r)\,dr  \ge \omega_{N-1} (2r_0)^{N-1} |T_k|.
$$
where $\omega_{N-1}$ denotes the euclidean surface measure of the unit sphere in $\R^N$. Thus (\ref{eq:tktozero}) holds, as claimed. From (\ref{eq:tktozero}) we deduce that
\begin{equation}
  \label{measure-estimate-3}
|R_k| \to r_0  \qquad \text{as $k \to \infty$.}
\end{equation}
Moreover,
\begin{equation}
  \label{eq:Srk-intersection}
S_{r,k} \cap \partial \Omega_k \not= \varnothing \qquad \text{for every $r \in R_k$}
\end{equation}
as a consequence of (\ref{eq:intersection-1}). Next we claim that
\begin{equation}
  \label{eq:measure-estimate-final}
\sigma_g(\Gamma_k) \ge |R_k| \qquad \text{for every $k \in \N$.}
\end{equation}
Here we shall need the assumption $N=2$. To derive (\ref{eq:measure-estimate-final}), we fix $k \in \N$, $\eps \in (0,\frac{r_0}{2})$ and set
$$
\Gamma_{k,\eps}:= \{x \in \partial \Omega_k\::\: 2 r_0 + \eps \le \dist(x,p_k) \le 3r_0 -\eps\} \subset \Gamma_k
$$
Since $\Omega_k$ is smooth and $\Gamma_{k,\eps}$ is compact,
 only finitely many (disjoint) path components $\Gamma_k^1,\dots,\Gamma_k^m$ of
$\Gamma_k$ intersect $\Gamma_{k,\eps}$. Let
$$
\beta_j^+:= \max \{\dist(p_k,q)\::\: q \in \Gamma_k^j \cap \Gamma_{k,\eps} \}\qquad \text{and}\qquad \beta_j^-:= \min \{\dist(p_k,q)\::\: q \in \Gamma_k^j \cap \Gamma_{k,\eps} \}
 $$
for $j=1,\dots,k$. By construction and (\ref{eq:Srk-intersection}) we then have
\begin{align*}
R_k \cap [2r_0 + \eps, 3r_0 -\eps] &\subset \{r \in [2r_0 + \eps, 3r_0 -\eps]\::\: \dist(x,p_k)=r \quad \text{for some $x \in
\Gamma_k^j \cap \Gamma_{k,\eps}$ and some $j$}\}\\
&\subset \{r \in [2r_0 + \eps, 3r_0 -\eps]\::\: \beta_j^- \le r \le  \beta_j^+ \; \text{for some $j$}\}
\end{align*}
and therefore
 \begin{equation}
  \label{eq:R_k-measure-estimate}
\bigl|R_k \cap [2r_0 + \eps, 3r_0 -\eps]\bigr| \le \sum_{j=1}^m (\beta_j^+- \beta_j^-).
 \end{equation}
Moreover, for every $j \in 1,\dots,m$ there is a smooth curve $\gamma_j:[0,1] \to \Gamma_k^j$ such that
$$
|\dot \gamma|_g>0\;\; \text{on $[0,1]$,}   \qquad \;\dist(p_k,\gamma(0))=\beta_j^- \;\qquad \text{and}\;\qquad \dist(p_k,\gamma(1))=\beta_j^+.
$$
Consequently,
\begin{align*}
\beta_j^+- \beta_j^- &= \int_0^1 \frac{d}{ds} \dist(p_k,\gamma(s))\,ds = - \int_0^1 \frac{1}{|\textrm{Exp}_{\gamma(s)}^{-1}(p_k)|_g}
\langle \textrm{Exp}_{\gamma(s)}^{-1}(p_k), \dot \gamma(s)\rangle_g \,ds\\
&\le \int_0^1 | \dot \gamma(s)|_g \,ds \le \sigma_g(\Gamma_k^j),
\end{align*}
the last inequality being a consequence of the fact that $\Gamma_k^j$ is a one-dimensional submanifold of $\M$. Here the assumption $N=2$ enters.
Combining this estimate with \eqref{eq:R_k-measure-estimate}, we deduce that
$$
\bigl|R_k \cap [2r_0 + \eps, 3r_0 -\eps]\bigr| \le \sum_{j=1}^m \sigma_g(\Gamma_k^j) \le \sigma_g(\Gamma_k).
$$
By considering the limit $\eps  \to 0$ we conclude that $|R_k| \le  \sigma_g(\Gamma_k)$, as claimed in (\ref{eq:measure-estimate-final}). Combining this inequality with (\ref{measure-estimate-3}) gives (\ref{measure-estimate-2}).
The proof is thus finished.
\end{pf}

We may now complete the proof of Theorem~\ref{sec:global-upper-bound-5} as follows (by contradiction):\\
By (\ref{eq-dim-2:11}), Lemma~\ref{sec:global-upper-bound-3}(ii) and the variational characterization of $\nu_2(\Omega_k,g)$ we have that
$$
\limsup_{k \to \infty} \nu_2(\Omega_k,g) \le \frac{\limsup \limits_{k \to \infty} \int_{\Omega_k} |\nabla f_k|_g^2 \,dv_g}{\liminf \limits_{k \to \infty} \int_{\partial \Omega_k} f_k^2\,d\sigma_g} \le
\frac{c_0^2 \lim \limits_{k \to \infty}  |\Omega_k|_g}{c_1} = 0,
$$
which contradicts (\ref{eq:dim2-asympt-1}). The proof of Theorem~\ref{sec:global-upper-bound-5} is thus finished.

\end{document}